\documentclass{article}
\usepackage{amsfonts,amsmath,amssymb,amsthm}
\usepackage{graphicx} 
\usepackage{tabulary}
\usepackage{todonotes}
\usepackage[margin=0.75in]{geometry}
\usepackage{biblatex}
\usepackage{algorithm}
\usepackage{algpseudocode}
\usepackage[T1]{fontenc}
\usepackage{authblk}

\usepackage{multicol}
\usepackage{todonotes}
\usepackage{booktabs}
\usepackage{hyperref}
\usepackage{algorithm}
\usepackage{algpseudocode}
\usepackage{tablefootnote}
\usepackage{subcaption}
\usepackage[capitalize]{cleveref}

\RequirePackage{tgtermes}
\RequirePackage{newtxtext}
\RequirePackage{newtxmath}
\RequirePackage{bm}
\RequirePackage{endnotes}

\usepackage{tikz}

\title{A Hybrid Decomposition Approach for Stochastic Unit Commitment with Combined-Cycle Generators}
\author[1*]{Rosemary Barrass}
\author[2]{Harsha Nagarajan}
\author[1]{Mathieu Tanneau}
\author[2]{Russell Bent}
\author[1]{Pascal Van Hentenryck}

\affil[1]{Georgia Institute of Technology}
\affil[2]{Los Alamos National Laboratory}
\affil[*]{Address correspondence to: rbarrass3@gatech.edu}

\bibliography{bib/AI4Opt-Rosemary}
\date{}

\begin{document}

\maketitle

\begin{abstract}%
The U.S. power grid is undergoing a paradigm shift as energy demand grows in scale and volatility.
In response to this growing need, the U.S. has increased adoption of combined-cycle generators (CCs).
CCs are fast-ramping generators that utilize variable configurations of combustion turbines (CTs) and steam turbines (STs) to achieve higher efficiency than traditional CTs.
For schedule optimization, modeling these CCs requires the addition of a large number of binary constraints and variables to Unit Commitment (UC) problem formulations.
This paper presents a novel hybrid Benders' (BD) and Dantzig-Wolfe (DW) decomposition algorithm, called CRG, for stochastic UC problems with CCs. 
CRG exploits the separability of the linear constraints in UC through BD and the integer CC constraints through DW.
A novel set of valid inequalities are proposed for significantly tightening the lower bound produced by CRG.
CRG is tested on the 935-generator FERC test data set, modified to include CC mode data.
Results demonstrate better primal solutions than BD on cases with at least 20 load scenarios.
CRG scales computationally better than Gurobi's branch-and-bound solver, which exceeds 64GB RAM allocations at 45 scenarios.
Results show that the proposed algorithm is a scalable approach for solving large-scale stochastic UC with CCs.
\end{abstract}%






\section{Introduction}
The U.S. power grid is in the midst of a major paradigm shift.
The rapid deployment of AI data centers and the increasing adoption of electric vehicles is fundamentally changing the distribution and scale of load across the country.
These major changes are leading to more volatility in load profiles and an increasing need for fast ramping generators.
These volatile load profiles can be modeled effectively in Unit Commitment (UC) models using stochastic formulations.
However, stochastic UC models are more computationally expensive than deterministic models and become increasingly intractable as more scenarios are added to the formulation.
Volatile systems like that of the modern grid require a large number of scenarios for the model to be robust to the worst cases.
System operators attempt to prepare for these worst cases by generating UC decisions every day under tight time constraints, making the trade-off between decision robustness and solve time critical to daily grid operations.
In addition to widely varying loads, system operators need to make decisions that will remain feasible for high ramping load scenarios.
To address these ramping needs, U.S. grid operators have adopted a high mix of combined-cycle generators (CCs).
As of 2023, CCs make up more than 20\% of the entire U.S. generation fleet (\cite{USElctrcGenByEngySrce2023,ElctrcPwrAnnual2023}).
These generators are fast ramping, generally low-cost natural gas generators that combine combustion turbines with steam turbines to make use of excess heat and improve operational efficiency.
These generators are also more complex to schedule than traditional combustion turbines because of their modularity: a CC can be configured to operate one or more turbines independently of each other.
From here on, the terms ``configuration'' and ``mode'' are used interchangeably to refer to a feasible operating configuration of the individual turbines of a CC.
The increased operating complexity of CCs is modeled by adding more binary variables and integer constraints into the standard formulation of UC.
It is well known that combinatorial optimization with generic solvers like Gurobi and CPLEX is computationally expensive, so it follows that UC formulations that include CC constraints and variables are more computationally expensive to solve.
Thus combining both a large, volatile stochastic UC formulation with CC modeling is a practical but challenging problem that is critical to modern grid operations.

\subsection{Literature Review}
While significant research has been devoted to solving both stochastic UC and deterministic UC with CCs, the integration of CCs into stochastic UC formulations remains understudied. 
This is largely due to the scarcity of publicly available data sets with detailed operational data on CCs, which limits the development and benchmarking of realistic models in academic settings.
Additionally, the combined challenges of a large number of binary constraints and variables in CC formulations plus the large number of load scenarios in stochastic formulations creates a problem that can be entirely intractable at real-world scale.
Therefore, to the best of the authors' knowledge, the current literature focuses on these topics completely separately.

While there is no literature on stochastic UC with CCs, deterministic UC with CCs is a well-researched topic with many proposed approaches.
A substantial portion of the literature on solving deterministic UC with CCs utilizes decomposition techniques, particularly Lagrangian Relaxation (LR) and Surrogate Lagrangian Relaxation (SLR), to manage computational complexity (\cite{AnEffcntApprchForUCAndEDwithCCandACPF2016}, \cite{ANvlDcmpstinAndCoordntnAprchForLrgDAUCwithCCs2018}, \cite{ShrtTermSchdlngOfCCs2004}, \cite{SLRAndBCforUCwithCCs2014}, \cite{PrceBsdUC:ACaseOfLRvsMIP2005}, \cite{AFstSltionMthdForLrgSclUCBsdOnLRAndDP2024}, \cite{ARLEmbddSLRMthdForFstSlvngUCPrblms2025}). 
These methods decompose the problem into tractable inner problems, enabling the solution of UC instances with many units over long planning horizons.
While LR and SLR are computationally more tractable than branch-and-cut (B\&C) algorithms, these methods do not guarantee convergence to the optimal solution but rather provide a solution to the approximation of the problem.

To address this optimality concern, researchers have also worked towards tighter formulations of UC with CCs to improve the computational efficiency of solving the full mixed integer linear program (MILP).
Many of these techniques involve clever solution space reduction techniques that exploit expert knowledge of real-world CC operation.
In (\cite{ASmplfdCCUnitMdlForMILPBasedUC2008}), the authors utilize a configuration-based model, which schedules physically feasible turbine combinations instead of individual turbines, and eliminate configurations that have cost curves which are always higher than all other configurations.
This approach, however, lacks the nuance and accuracy of individual component modeling.
In (\cite{CompModeModforShortTermSchedofCCs2009}), the authors present a component-based model that takes advantage of the more accurate fuel-generation curve modeling to add additional power augmenting features, like duct burners, to improve the UC model accuracy.
The results show that this modeling approach is tractable for the IEEE 118-bus system.
However, this component-based scheduling cannot be used directly in energy markets for day-ahead bidding.
In (\cite{HybridCmpntAndConfigModelForCCsInUC2018}), the authors address this mismatch of scheduling information by presenting a hybrid component-configuration model which uses configuration-based pricing curves with component-based physical constraints.
While this model is ideal for accuracy and efficacy in the real-world, this formulation takes longer to solve on the IEEE-118 bus system than the configuration-based model.
Finally, (\cite{UCwithFlxbleGenUnits2005}) uses a combination of BD and LR to exploit the separability of the binary and linear constraints in scheduling flexible fuel CCs.
This model is an aggregate UC model, with all CCs modeled as single generating units with no modes, but there is a similar level of additional complexity added to the system in the form of fuel operating modes.
Across all of these formulations and decomposition methods, the balance between computational efficiency and model accuracy is key, with the main setbacks coming from the large number of binary variables and constraints that each CC requires.

Additionally, nearly all of the literature is tested using a relatively small system, the modified IEEE-118 bus system, with the exception of a few authors with access to proprietary data from real system operators.
While the real-world cases are impressive, limited data transparency hinders reproducibility and rigorous benchmarking.
Consequently, solving UC problems with CCs remains a challenging and computationally intensive task, highlighting the need for further research and better data availability in this area.

We draw on ideas beyond the UC literature to leverage a key property of UC with CCs: per-generator separability. 
The added binaries required for accurate CC modeling increase complexity, making the problem well-suited to column generation, as in job-shop scheduling and vehicle routing. 
Similar approaches have succeeded for stochastic UC without CCs (\cite{StochUCProblem2004}). 
More broadly, pairing column generation with constraint/row generation is an empirically observed scalability method for MILPs with block-diagonal structures (\cite{BDandCRGforSlvngLrgSclLPsWithClDepRows2018,SlvngTwoStgeROPrblmsUsingCCGMthd2013,CRGMthdToABtchMachneSchedPrblm2010,CmbngDWandBDtoSlveALrgSclNclrOutgePlnningPrblm2022, EfctvHybrdDecmpApprchToSlvNtwrkCnstrndStochUC2024}).

\subsection{Contributions and Outline}
This paper introduces a hybrid column-row generation (CRG) framework designed to tackle the computational challenges posed by both large sets of load scenarios and the large number of binary variables and constraints in UC with CC.
The decomposition method combines Benders' Decomposition (BD), which is often used to improve solve times in stochastic programming, and Dantzig-Wolfe Decomposition (DW), which is used in column generation methods for reducing the size of integer problem search spaces.
This methodology is applied to a stochastic UC (SUC) formulation with CCs and varying load scenarios; an increasingly relevant, critical real-world use-case that has not yet been explored in literature as a single problem.
The proposed combination of decomposition techniques is a novel approach to solving a previously computationally intractable size of SUC by utilizing the strongest aspects of both decomposition techniques.
The framework is further accelerated by a novel family of valid inequalities that strengthen the lower bounds of scenario recourse costs.
This paper demonstrates the efficacy of this hybrid decomposition on large-scale, open-source power-system benchmarks, outperforming both a traditional stochastic programming technique and a commercial branch-and-bound solver.

The paper is organized as follows: Section \ref{sec:methodology} describes the SUC formulation, decomposition techniques, and acceleration techniques for the proposed methodology, Section \ref{sec:computational-results} compares the results of testing the proposed algorithm, BD, LR, and a state-of-the-art solver on large-scale test cases, and Section \ref{sec:conclusion} presents the concluding remarks and future directions for this work.

\section{Methodology}\label{sec:methodology}

Motivated by decades of progress in column and row generation for MILPs, this paper introduces a novel CRG framework for SUC with CCs, and demonstrates its efficacy on large-scale, open-source power-system benchmarks.

Methodologically, CRG is a two-level decomposition. 
A BD splits SUC into an outer problem for generator commitment and a linear inner problem for networked dispatch, reserves, and scenario/contingency recourse. 
Within the outer problem, a DW reformulation disaggregates generator configuration constraints (minimum up/down, ramping, startup/shutdown) into smaller-scale per-generator pricing inner problems that generate columns on demand, while row generation adds violated inner problem constraints. 
The proposed algorithm exploits the separability of the problem for faster solution times.


\subsection{Full UC with CC Formulation}
\newcommand{\cf}{\mathbf{u}}
\newcommand{\su}{\mathbf{v}}
\newcommand{\sd}{\mathbf{w}}
\newcommand{\pg}{\mathbf{p}}
The UC problem is a generator scheduling problem which power grid operators solve every day to maintain grid reliability.
The formulation used in this paper is a standard UC problem with mode-based CC modeling, which means that it includes engineering constraints for the feasible operation of each generator and mode.
Throughout this paper, variables are denoted with bold characters.
The standard model is made stochastic by the inclusion of multiple load scenarios.
For ease of reading, the load scenario indices are omitted; the paper considers a two-stage stochastic formulation where first-stage variables comprise (binary) commitment, startup and shutdown variables, and second-stage variables comprise (continuous) energy dispatch, reserve, load-imbalance, and reserve-shortfall variables.

Let $\mathcal{T}$ and $\mathcal{G}$ denote the set of all time steps and of all generators, respectively.
For generator $g \in \mathcal{G}$, let $\mathcal{M}_{g}$ denote this generator's set of modes.
For each mode $m \in \mathcal{M}_g$, let $\mathcal{L}_m$ denote the set of piecewise energy-dispatch points used to approximate mode $m$'s nonlinear dispatch curve, and let $p^l_{g,m,t}$ denote the total output at point $l\in\mathcal{L}_m$.
Next, let binary variables $\mathbf{u}_{g, m, t}, \mathbf{v}_{g,m,t}, \mathbf{w}_{g,m,t} \in \{0, 1\}$ denote the commitment, startup and shutdown status of generator $g$ and mode $m$ at time $t$.
Continuous variables $\mathbf{p}_{g, m, t}$ and $\mathbf{r}_{g, m, t}$ model the marginal energy dispatch and spinning reserve quantity, respectively.
Let the marginal energy dispatch be defined as the incremental dispatch above the minimum limit $\underline{p}_{g, m, t}$.
Let the spinning reserve provided by generator $g$ in mode $m$ at time $t$ be defined as the available headroom that generator $g$ can provide in one time step.
$\mathbf{r}_{g, m, t}$ is only provided by thermal generators and not by wind, solar, or storage units.
The continuous variables $\boldsymbol{\lambda}_{g, m, l, t}$ represent the portion of energy dispatch from piecewise linear dispatch approximation point $l$ for generator $g$ and mode $m$ at time $t$.
Generator minimum and maximum dispatch limits are denoted by $\underline{p}_{g, m, t}$ and $\overline{p}_{g, m, t}$, respectively.
Similarly, generator ramping limits are denoted by $\overline{r}^{\text{up}}_{g, m, t}, \overline{r}^{\text{down}}_{g, m, t}, \overline{r}^{\text{startup}}_{g, m, t}$ and $\overline{r}^{\text{shutdown}}_{g, m, t}$, which capture upward, downward, startup and shutdown ramping limits, respectively.
All variables and parameters can be represented at the unit-level by removing the index $m$, e.g. $\cf_{g, t}$ is the unit-level commitment status of generator $g$ at time $t$.
Finally, the nonnegative variables $\boldsymbol{\delta}^{+}_{t}$,
$\boldsymbol{\delta}^{-}_{t}$, and $\boldsymbol{\delta}^{R}_{t}$ denote load
shortage, excess generation, and reserve shortfall at time $t$, respectively.

The objective function minimizes generation and penalty costs, expressed as
\begin{align}
    \label{eq:UC:obj_function}\tag{\text{Full UC}}
    \min_{\cf, \pg} & \quad \sum_{g \in \mathcal{G}, m \in \mathcal{M}_g, t \in \mathcal{T}} \sum_{l \in \mathcal{L}_m} (c^{\text{dispatch}}_{g, m, l, t} - c^{\text{dispatch}}_{g, m, 1, t}) * \boldsymbol{\lambda}_{g, m, l, t}\\  
    \quad & + \sum_{g \in \mathcal{G}, m \in \mathcal{M}_g, t \in \mathcal{T}} c^{\text{dispatch}}_{g, m, 1, t} * p^1_{g,m,t} * \cf_{g, m, t} \nonumber\\ 
    \quad & + \sum_{t \in \mathcal{T}} c^{\text{penalty}}_t
    * \left(\boldsymbol{\delta}^{+}_{t}+\boldsymbol{\delta}^{-}_{t}
    +\boldsymbol{\delta}^{R}_{t}\right) \nonumber
\end{align}
where $c^{\text{dispatch}}_{g, m, l, t}$ is the $l$-th piecewise linear approximation of the cost of energy dispatch for generator $g$ in mode $m$ at time $t$ and $c^{\text{penalty}}_{t}$ is the common penalty cost for load shortage, excess generation, and reserve shortfall at time $t$.

The startup/shutdown logic constraints, minimum up time, and minimum down time constraints for each generator, mode, and time step read
\begin{align}
    \label{eq:UC:commitment}
    \cf_{g,m,t} &= \cf_{g,m,t-1} + \su_{g,m,t} - \sd_{g,m,t} \quad \forall g, m, t \in \mathcal{G}, \mathcal{M}_g, \mathcal{T}\\
    \label{eq:UC:min_run_time}
    \cf_{g,m,t} &\geqslant \sum_{i=t - T^{minup}_{g,m}+1}^{t} \su_{g,m,i} \quad \forall g, m, t \in \mathcal{G}, \mathcal{M}_g, \mathcal{T}\\
    \label{eq:UC:min_dwn_time}
    \cf_{g,m,t} &\leqslant 1 - \sum_{i=t - T^{mindn}_{g,m}+1}^{t} \sd_{g,m,i} \quad \forall g, m, t \in \mathcal{G}, \mathcal{M}_g, \mathcal{T}
\end{align}
where $T^{minup}_{g,m}, T^{mindn}_{g,m}$ are the minimum up and minimum down time, respectively, for each generator $g$ and mode $m$.
In \eqref{eq:UC:commitment}, the generator status is turned on or off based on whether or not a startup or shutdown has occurred.
If neither have occurred, then consecutive decisions must be equal.
In \eqref{eq:UC:min_run_time} and \eqref{eq:UC:min_dwn_time}, the startups and shutdowns are limited to only occur once within the minimum up time and minimum town time range, respectively.

Dependent mode startup and commitment is restricted by its supporting mode commitment and startup status, respectively.
\begin{align}
    \label{eq:UC:dep-mode-startup}
    \mathbf{v}_{g, m_d, t} & \leqslant \mathbf{u}_{g, m_s, t} \quad \forall g, (m_d, m_s), t \in \mathcal{G}, \mathcal{M}_g, \mathcal{T}\\
    \label{eq:UC:dep-mode-logic}
    \mathbf{u}_{g, m_d, t} & \leqslant 1 - \mathbf{v}_{g, m_s, t} \quad \forall g, (m_d, m_s), t \in \mathcal{G}, \mathcal{M}_g, \mathcal{T}
\end{align}
In \eqref{eq:UC:dep-mode-startup} and \eqref{eq:UC:dep-mode-logic}, the mode subscript ($d$) denotes a dependent mode and ($s$) a supporting mode. 
For CCs, a base mode is any configuration into which the unit can start that requires no supporting modes to be online.
A dependent mode is any mode that requires another mode (i.e. its supporting mode) to be online simultaneously.
The power generated by the modes is additive, so generator $g$ has total
generation $\sum_{m\in\mathcal{M}_g}
(\underline{p}_{g,m,t}\cf_{g,m,t}+\mathbf{p}_{g,m,t})$ at time $t$.

For each generator, at most one base mode may be online at any time, yielding
\begin{align}
    \label{eq:UC:single-mode-logic}
    \sum_{m_{\text{base}} \in \mathcal{M}_g} \mathbf{u}_{g, m_{\text{base}}, t} &  =\mathbf{u}_{g, t} \quad \forall g, t \in \mathcal{G}, \mathcal{T},
\end{align}
where the subscript \textit{base} denotes any mode that is a supporting mode but never a dependent mode.

The ramp-up and ramp-down limits are given by:
\begin{align}
    \label{eq:UC:ramp-up-limits}
    \mathbf{p}_{g, m ,t} + \mathbf{r}_{g, m, t} - \mathbf{p}_{g, m, t-1} \leqslant \ \overline{r}^{\text{up}}_{g, m, t} \quad \forall g, m, t \in \mathcal{G}, \mathcal{M}_g, \mathcal{T}\\
    \label{eq:UC:ramp-dn-limits}
    \mathbf{p}_{g, m, t-1} - \mathbf{p}_{g, m ,t} \leqslant \ \overline{r}^{\text{down}}_{g, m, t} \quad \forall g, m, t \in \mathcal{G}, \mathcal{M}_g, \mathcal{T}
\end{align}
These constraints ensure that dispatch levels do not change too rapidly between time steps for physical generator limitations.

The startup and shutdown ramping limits are given by:
 \begin{align}
    \mathbf{p}_{g, m, t} + \mathbf{r}_{g, m, t} - (\overline{p}_{g, m, t} - \underline{p}_{g, m, t})\mathbf{u}_{g, m, t} \quad  & \quad \nonumber\\
    + \max \left\{(\overline{p}_{g, m ,t} - \overline{r}^{\text{startup}}_{g, m, t}), 0\right\}\mathbf{v}_{g, m, t} \leqslant ~ 0 & \quad \forall g, m, t \in \mathcal{G}, \mathcal{M}_g, \mathcal{T} \label{eq:UC:su-ramp-limits}\\
    \mathbf{p}_{g, m, t} + \mathbf{r}_{g, m, t} - (\overline{p}_{g, m, t} - \underline{p}_{g, m, t})\mathbf{u}_{g, m, t}  \quad & \quad \nonumber \\
    + \max \left\{(\overline{p}_{g, m ,t} - \overline{r}^{\text{shutdown}}_{g, m, t}), 0\right\}\mathbf{w}_{g, m, t+1} \leqslant ~ 0 & \quad  \forall g, m, t \in \mathcal{G}, \mathcal{M}_g, \mathcal{T} \label{eq:UC:sd-ramp-limits}
 \end{align}
 These constraints ensure the change in generation from online to offline, and vice versa, does not exceed the physical limitations of the generators.

The system-wide generation must meet system power demand:
\begin{align}
    \label{eq:UC:power-balance}
    \sum_{g,m \in \mathcal{G}, \mathcal{M}_g}
    \left(\cf_{g,m,t}\underline{p}_{g,m,t}+\mathbf{p}_{g,m,t}\right)
    +\boldsymbol{\delta}^{+}_{t}-\boldsymbol{\delta}^{-}_{t}
    =d_t \quad \forall t\in\mathcal{T}
\end{align}
where $\boldsymbol{\delta}^{+}_{t}$ and $\boldsymbol{\delta}^{-}_{t}$ account
for load shortage and excess generation, respectively.
The system-wide reserve requirement must also be met:
\begin{align}
    \label{eq:UC:reserve-requirement}
    \sum_{g,m \in \mathcal{G}, \mathcal{M}_g}\mathbf{r}_{g,m,t}
    +\boldsymbol{\delta}^{R}_{t}\geqslant R_t
    \quad \forall t\in\mathcal{T}
\end{align}
where $R_t$ is the system-wide reserve requirement at time $t$ and
$\boldsymbol{\delta}^{R}_{t}$ accounts for reserve shortfall and provides complete recourse.

Finally, nonlinear energy dispatch curves are approximated using piecewise linear approximations.
The dispatch and commitment relations for the piecewise weights are
\begin{align}
    \label{eq:UC:piecewise-dispatch}
    \mathbf{p}_{g,m,t}
    &=\sum_{l\in\mathcal{L}_m}
      \left(p^l_{g,m,t}-\underline{p}_{g,m,t}\right)
      \boldsymbol{\lambda}_{g,m,l,t}
    && \forall g,m,t\in\mathcal{G},\mathcal{M}_g,\mathcal{T},\\
    \label{eq:UC:piecewise-approximation}
    \sum_{l\in\mathcal{L}_m}\boldsymbol{\lambda}_{g,m,l,t}
    &=\cf_{g,m,t}
    && \forall g,m,t\in\mathcal{G},\mathcal{M}_g,\mathcal{T},\\
    \mathbf{p}_{g,m,t},\ \mathbf{r}_{g,m,t}, \ \boldsymbol{\lambda}_{g,m,l,t},&  \ \boldsymbol{\delta}^{+}_{t},\ \boldsymbol{\delta}^{-}_{t}, \
    \boldsymbol{\delta}^{R}_{t}
    \geqslant 0
    && \forall g,m,l,t\in
       \mathcal{G},\mathcal{M}_g,\mathcal{L}_m,\mathcal{T}.
\end{align}

\subsection{General formulation}

\newcommand{\x}{\mathbf{x}}
\newcommand{\y}{\mathbf{y}}

For clarity, we present the SUC model in general form:
\begin{subequations}
\begin{align}
    \min_{\x, \y} \quad 
        & \label{eq:UC:MIP:general:obj}
         c^{\top} \x + q^{\top} \y \tag{\text{General SUC}}\\
    \text{s.t.} \quad
        & \label{eq:UC:MIP:general:combinatorial}
        \x \in \mathcal{X}\\
        & \label{eq:UC:MIP:general:mixed}
        T \x + W \y \geqslant h 
\end{align}
\end{subequations}
where constraint \eqref{eq:UC:MIP:general:combinatorial} captures generators' scheduling constraints (e.g. \eqref{eq:UC:commitment}-\eqref{eq:UC:single-mode-logic}) and constraint \eqref{eq:UC:MIP:general:mixed} captures system-wide linking constraints (e.g. \eqref{eq:UC:power-balance}-\eqref{eq:UC:reserve-requirement}).
The following methodology is presented in this general form to illustrate that it can be expanded to include line flow constraints, additional reserve types, etc.
Additionally, the scenario indices $s$ for the system-wide constraints are omitted for simplicity.
Any additional constraints not included in \eqref{eq:UC:obj_function} can be sorted into the categories described above and decomposed in the same manner while maintaining the feasibility of this methodology.
The simple SUC formulation described in this paper is used for the sake of readability.

\subsection{Benders' Decomposition}

Separating first- and second-stage variables $\x$ and $\y$ yields the Benders' outer problem
\begin{subequations}
\begin{align}
    \min_{\x, \boldsymbol{\eta}} \quad
        & c^{\top} \x + \boldsymbol{\eta} \label{eq:UC:BD:general:master}\\
    \text{s.t.} \quad
        & \x \in \mathcal{X}\label{eq:UC:BD:general:x-polyhedron}\\
        & \boldsymbol{\eta} \geqslant Q(\x)
            \label{eq:UC:BD:general:value_function}
\end{align}
\end{subequations}
where $Q(\x)$ is the second-stage value function, defined as
\vspace{-0.5cm}
\begin{subequations}
\label{eq:UC:BD:general:subproblem}
\begin{align}
    Q(\x) = 
    \min_{\y} \quad
        & q^{\top} \y \\
    \text{s.t.} \quad
        & W \y \geqslant h - T \x
\end{align}
\end{subequations}
Each second-stage function represents a different net load scenario, though the indices are omitted for simplicity.
The nonlinear constraint \eqref{eq:UC:BD:general:value_function} is lower-approximated in the outer problem using Benders' cuts of the form
$
    \boldsymbol{\eta} \geqslant \gamma^{\top}(h - T \x),
$
where $\gamma \in \Gamma$ and $\Gamma = \{\gamma \geqslant 0 \, | \, W^{\top} \gamma = q \}$ is the set of all dual-feasible solutions for the BP-IP.
This gives a final outer problem formulation
\begin{subequations}
\begin{align}
    \min_{\x, \boldsymbol{\eta}} \quad
        & c^{\top} \x + \boldsymbol{\eta} \label{eq:UC:BD:general:master-with-cuts}\\
    \text{s.t.} \quad
        & \x \in \mathcal{X}\\
        & \boldsymbol{\eta} \geqslant \gamma^{\top}(h - T \x),
        && \forall \gamma \in \bar{\Gamma}
            \label{eq:UC:BD:general:cuts}
\end{align}
\end{subequations}
where $\bar{\Gamma} \subseteq \Gamma$ is identified with the set of cuts currently considered in the outer problem.
Note that the paper's implementation uses a so-called \emph{multi-cut} formulation, where cuts are added separately for each scenario.
This approach is known to improve convergence.

\subsection{Dantzig-Wolfe Decomposition}

The Dantzig-Wolfe (DW) decomposition for generator's scheduling constraints is obtained by convexifying constraints \eqref{eq:UC:MIP:general:combinatorial}.
Namely, letting $\Omega$ denote the set of extreme points of $\mathcal{X}$, constraint \eqref{eq:UC:MIP:general:combinatorial} is replaced with 
$
    \x = \sum_{\omega \in \Omega} \boldsymbol{\mu}_{\omega} \omega,
$
where $\mu$ are non-negative weights that sum to $1$.
This yields the DW outer problem
\begin{subequations}
\label{eq:UC:DW:general:master}
\begin{align}
    \text{OP}_{DW}
    \quad
    \min_{\x, \boldsymbol{\mu}, \boldsymbol{\eta}} \quad
    & c^{\top} \x + \boldsymbol{\eta} \\
    \text{s.t.} \quad
    & \sum_{\omega \in \Omega} \boldsymbol{\mu}_{\omega} \omega - \x = 0 
        && [\pi] \label{eq:UC:DW:general:linear-combination-constraint}\\
    & \sum_{\omega \in \Omega} \boldsymbol{\mu}_{\omega} = 1 
        && [\sigma]
        \label{eq:UC:DW:general:linear-coefficient-sum-constraint}\\
    & \boldsymbol{\mu} \geqslant 0\\
    & \boldsymbol{\eta} \geqslant \gamma^{\top}(h - T \x)
        && \forall \gamma \in \bar{\Gamma}
        \label{eq:UC:general:BendersCut}
\end{align}
\end{subequations}
where $\pi$ and $\sigma$ are the duals of constraints \eqref{eq:UC:DW:general:linear-combination-constraint} and \eqref{eq:UC:DW:general:linear-coefficient-sum-constraint}, respectively.
Note that, in order to formulate Benders' cuts efficiently, the present formulation keeps variables $\x$ in the DW outer problem.

New columns are identified by solving the pricing inner problem
\begin{align}
    \label{eq:UC:DW:general:subproblem}
    \text{IP}_{DW}(\pi, \sigma):
    \quad
    \min_{\x \in \mathcal{X}} \quad
        - \pi^{\top}\x - \sigma
\end{align}
This inner problem can then be separated into $\mathcal{G}$ independent problems, one per generator.
Thus, the final decomposition of the problem is separated into 3 parts: a DW outer problem, per-generator pricing inner problems, and per-load-scenario inner problems.

\subsection{Column-Row Generation Algorithm}
To solve this series of outer and inner problems, a CRG algorithm is employed (see Appendix \cref{alg:CRG}).
This algorithm iteratively adds columns from the Pricing IP and cuts from the Benders' IP.
An initial set of columns is generated by seeding the DW inner pricing problems with the dual solution of the linear relaxation  \eqref{eq:UC:obj_function}.
The algorithm proceeds in a loop until there are no more columns or rows to generate, or the desired optimality gap is achieved.
At each iteration, valid upper and lower bounds are computed as follows.
A valid upper bound is obtained as
\begin{align}
    \label{eq:CRG:upper-bound}
    \text{UB}_{CRG}(\omega^*, \y) \quad = c^{\top}\omega^* + q^{\top}\y
\end{align}
where $\omega^*$ is the best schedule in the current set of schedules.
The lower bound is updated with the duals of the DW inner and outer problems:
\begin{align}\label{eq:CRG:lower-bound}
    \text{LB}_{CRG}(\x, \pi, \sigma, \text{DualObj}\eqref{eq:UC:DW:general:master}) \quad
    = -\pi^{\top}\x - \sigma + \text{DualObj}\eqref{eq:UC:DW:general:master}
\end{align}
where $\text{DualObj}\eqref{eq:UC:DW:general:master}$ is the objective value of the dual problem of \eqref{eq:UC:DW:general:master}.
The optimality gap is calculated by subtracting the lower bound from the upper bound and dividing this by the upper bound.
Once one of these criteria is met, the outer problem is converted back into an integer program with the generated columns and cuts and the final integer solution is found.

\subsection{Stabilization and Acceleration Techniques}
\label{subsec:stabilization}
Without stabilization, BD-based algorithms can suffer from slow convergence due to the somewhat erratic nature of the primal solution.
This behavior can be mitigated using in-out separation, as described in \cite{ImplmntngAutoBDInAMdrnMIPSlvr2020}.
For this paper, an in-out separation technique is applied to both BD and CRG to stabilize the primal solution and improve convergence.

\begin{figure}
    \centering
    \includegraphics[width=0.6\linewidth]{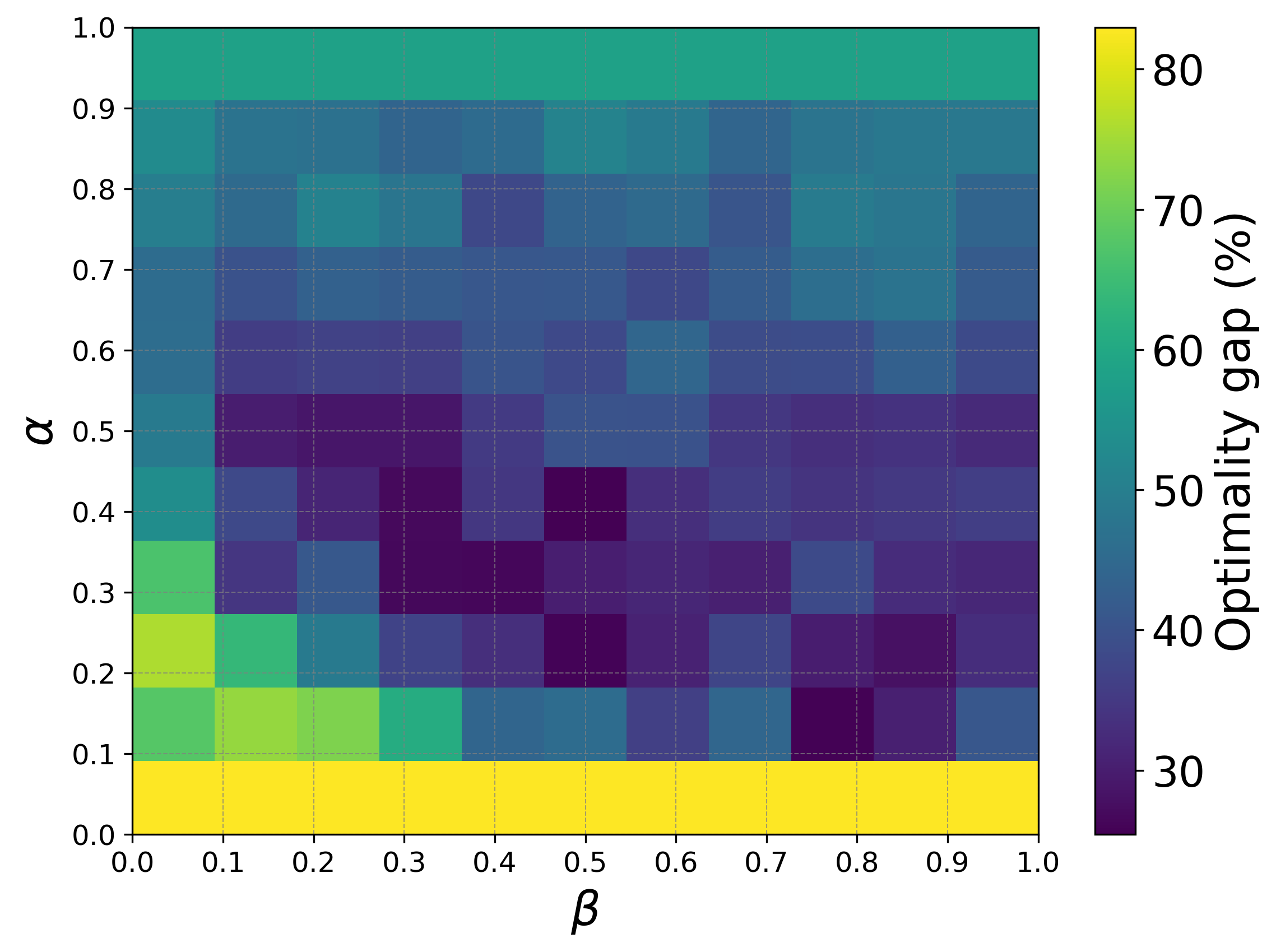}
    \caption{Sensitivity analysis on the 30-load scenario case terminated after 120s using varying $\alpha$ and $\beta$ values. The heatmap shows the optimality gap (\%) for each $(\alpha,\beta)$ pair, with grid lines (grey) indicating the evaluated combinations. Setting $\alpha$ and $\beta$ to $1$ is equivalent to not using in–out separation in Section \ref{subsec:stabilization}.}
    \label{fig:sensitivity-analysis}
\end{figure}

In-out separation uses a feasible point in the interior of the feasible region of the OP to provide a better solution to the inner problems than just using the previous outer problem solution.
In general, instead of solving the second stage of BD with the outer problem solution $x^{*k}$ and previous core point $\bar{x}^{k-1}$ as input and a variable $\alpha$ is chosen such that $0 \leqslant \alpha \leqslant 1$ to calculate the following:
\begin{equation}\label{eq:simple-io-separation}
    x_{io}^k = \alpha x^{*k} + (1-\alpha)\bar{x}^{k-1}
\end{equation}
The inner problems are then solved with the input $x_{io}^k$.
This technique is sensitive to the choices of $\alpha$ and $\beta$, as can be seen in \cref{fig:sensitivity-analysis}.
This sensitivity analysis was performed on a PGLIB-UC instance with 30 load scenarios, 934 thermal generators, and 12 time steps. Each run was terminated after 120 s. Optimality gaps between the best incumbent solution and the lower bound are reported in percent, computed as the difference divided by the best solution. Gaps are smallest (down to $\approx$25\%) for moderate values of $\alpha$ and $\beta$ (around 0.3–0.5), and much larger ($\approx$58–83\%) at the extremes $(\alpha=0)$ and $(\alpha=\beta=1)$.
So, moderate in–out separation clearly gives better solution quality than using no separation.

From this analysis, the best values for $\alpha$ and $\beta$ were chosen as 0.4 and 0.5, respectively.

Column generation (CG) algorithms can also suffer from adverse tailing-off effects due to erratic dual behavior.
To combat this, it is well-known that applying a box constraint to the duals can improve the stability of CG algorithms (\cite{BrnchAndPrce2026}).
For this paper, a basic dual box constraint is applied to the duals of each constraint in the OP before solving the dual problem to obtain the optimal duals for generating schedules.
This mitigates some of the tailing off effect that adversely impacts the convergence of the CG aspects of the CRG algorithm.
It is important to note that the size of these dual box constraints is problem-dependent.

For the CRG and BD tests, this approach was accelerated by warm-starting with a core point generated by solving the deterministic UC for the scenario with the highest single-hour load.
This in-out separation point is used to solve the inner problems and the core point is updated.
The core point is updated using another parameter $\beta$ such that $0 \leqslant \beta \leqslant 1$:
\begin{equation}\label{eq:corepoint}
    \bar{x}^k = \beta x_{io}^k + (1-\beta)\bar{x}^{k-1}
\end{equation}
This way, the in-out separation point can be chosen more aggressively, even when the outer solution is approaching optimal.
The branch-and-bound (B\&B) solver is also warm-started with the deterministic highest single-hour load scenario solution to provide a fair comparison.

Along with necessary stabilization techniques, CRG takes advantage of the large number of pricing IPs by dynamically adding a percentage of the negative reduced cost columns at each iteration.
This heuristic technique can reduce the number of iterations necessary for algorithm convergence at the cost of potentially overwhelming the OP.
For this reason, as with the stabilization techniques, a balance must be struck when choosing the correct number of columns to add at each iteration.

\subsection{Valid Inequality-based Acceleration}
\label{subsec:hybrid-valid-inequalities}
We further strengthen the multi-cut DW outer problem with scenario- and
time-specific lower bounds on second-stage dispatch and penalty costs. Let
$\mathcal{S}$ denote the set of load scenarios, and let
$\mathcal{G}^{R}\subseteq\mathcal{G}$ denote the reserve-capable generators.
The bounds exploit the copper-plate, lossless balance in
\eqref{eq:UC:power-balance}, nonnegative load-imbalance and reserve-shortfall
variables charged at the common unit penalty $c_t^{\text{penalty}}$, and
balance prices no larger than that penalty. Scenario indices on second-stage
variables such as $\boldsymbol{\lambda}$ are omitted for readability.
For $l\in\mathcal{L}_m$, let $p^{l}_{g,m,t}$ be the total output at piecewise
point $l$, so that the mode output is
$\sum_l p^{l}_{g,m,t}\boldsymbol{\lambda}_{g,m,l,t}$, and define
\begin{align}
    \Delta c^{l}_{g,m,t}
    =c^{\text{dispatch}}_{g,m,l,t}
    -c^{\text{dispatch}}_{g,m,1,t}.
\end{align}
Because $\boldsymbol{\lambda}_{g,m,l,t}\geqslant 0$ and
$\sum_{l\in\mathcal{L}_m}\boldsymbol{\lambda}_{g,m,l,t}
=\cf_{g,m,t}$ by \eqref{eq:UC:piecewise-approximation}, the incremental
dispatch cost is bounded below by
\begin{align}
    \Phi_t(\cf)
    &=\sum_{g\in\mathcal{G}}\sum_{m\in\mathcal{M}_g}
      b_{g,m,t}\cf_{g,m,t}, \quad 
    b_{g,m,t}
    =\min_{l\in\mathcal{L}_m}\left\{\Delta c^{l}_{g,m,t}\right\}.
    \label{eq:CRG:dispatch-cost-floor}
\end{align}
For nondecreasing cost curves, $b_{g,m,t}=0$ because the first piecewise
point is included in the minimum. The more general definition remains valid
for nonmonotone input data.

For each scenario and period, introduce an unrestricted lower-bound variable
$\boldsymbol{\xi}_{s,t}$ and impose
\begin{align}
    \boldsymbol{\eta}_s
    &\geqslant \sum_{t\in\mathcal{T}}\boldsymbol{\xi}_{s,t},
    && \forall s\in\mathcal{S},
    \label{eq:CRG:recourse-cost-aggregation}\\
    \boldsymbol{\xi}_{s,t}
    &\geqslant \Phi_t(\cf),
    && \forall s\in\mathcal{S},\ t\in\mathcal{T}.
    \label{eq:CRG:dispatch-floor-inequality}
\end{align}

\paragraph{Price-based bound.}
Let $\widehat{\varphi}_t$ be the power-balance dual obtained from the LP
relaxation for the load scenario attaining the largest single-period demand,
projected onto
$0\leqslant\widehat{\varphi}_t\leqslant c_t^{\text{penalty}}$. Define
\begin{align}
    \kappa_{g,m,t}(\widehat{\varphi}_t)
    =\min_{l\in\mathcal{L}_m}
    \left\{
        \Delta c^{l}_{g,m,t}
        -\widehat{\varphi}_t p^{l}_{g,m,t}
    \right\}.
\end{align}
Then the following inequality is valid for every scenario and period:
\begin{align}
    \boldsymbol{\xi}_{s,t}
    \geqslant
    \widehat{\varphi}_t d_{s,t}
    +\sum_{g\in\mathcal{G}}\sum_{m\in\mathcal{M}_g}
      \kappa_{g,m,t}(\widehat{\varphi}_t)\cf_{g,m,t},
    \quad \forall s\in\mathcal{S},\ t\in\mathcal{T}.
    \label{eq:CRG:price-bound}
\end{align}
Validity requires only the displayed price interval; the LP dual is used to
select an informative price. Pricing \eqref{eq:UC:power-balance} by
$\widehat{\varphi}_t$ and using the convex-combination relation above bounds
each committed mode by $\kappa_{g,m,t}(\widehat{\varphi}_t)$. The residual
coefficients of $\boldsymbol{\delta}^{+}_{t}$ and
$\boldsymbol{\delta}^{-}_{t}$ are, respectively,
$c_t^{\text{penalty}}-\widehat{\varphi}_t$ and
$c_t^{\text{penalty}}+\widehat{\varphi}_t$, and are therefore nonnegative.
Unlike a capacity-only bound,
\eqref{eq:CRG:price-bound} can remain positive for capacity-adequate
commitments.

\paragraph{Physical bounds.}
Let minimum output, maximum capacity, and reserve headroom be
\begin{subequations}
\begin{align}
    P_t^{\min}(\cf)
    &=\sum_{g\in\mathcal{G}}\sum_{m\in\mathcal{M}_g}
      \underline{p}_{g,m,t}\cf_{g,m,t},\\
    P_t^{\max}(\cf)
    &=\sum_{g\in\mathcal{G}}\sum_{m\in\mathcal{M}_g}
      \overline{p}_{g,m,t}\cf_{g,m,t},\\
    H_t(\cf)
    &=\sum_{g\in\mathcal{G}^{R}}\sum_{m\in\mathcal{M}_g}
      (\overline{p}_{g,m,t}-\underline{p}_{g,m,t})\cf_{g,m,t}.
\end{align}
\end{subequations}
For $g\in\mathcal{G}^{R}$, let
\begin{align}
    a^{\mathrm{SU}}_{g,m,t}
    &=\max\left\{
        \overline{p}_{g,m,t}-\overline{r}^{\text{startup}}_{g,m,t},0
      \right\}, \quad 
    a^{\mathrm{SD}}_{g,m,t}
    =\max\left\{
        \overline{p}_{g,m,t}-\overline{r}^{\text{shutdown}}_{g,m,t},0
      \right\},
\end{align}
and define the startup- and shutdown-limited capacities
\begin{subequations}
\begin{align}
    P_t^{\mathrm{SU}}(\cf,\su)
    &=P_t^{\max}(\cf)
      -\sum_{g\in\mathcal{G}^{R}}\sum_{m\in\mathcal{M}_g}
       a^{\mathrm{SU}}_{g,m,t}\su_{g,m,t},\\
    P_t^{\mathrm{SD}}(\cf,\sd)
    &=P_t^{\max}(\cf)
      -\sum_{g\in\mathcal{G}^{R}}\sum_{m\in\mathcal{M}_g}
       a^{\mathrm{SD}}_{g,m,t}\sd_{g,m,t+1},
      && t<|\mathcal{T}|.
\end{align}
\end{subequations}
The physical recourse inequalities are
\begin{subequations}
\label{eq:CRG:physical-recourse-bounds}
\begin{align}
    \boldsymbol{\xi}_{s,t}
    &\geqslant \Phi_t(\cf)
      +c_t^{\text{penalty}}
       \left(d_{s,t}+R_t-P_t^{\mathrm{SU}}(\cf,\su)\right),\\
    \boldsymbol{\xi}_{s,t}
    &\geqslant \Phi_t(\cf)
      +c_t^{\text{penalty}}
       \left(d_{s,t}+R_t-P_t^{\mathrm{SD}}(\cf,\sd)\right),
      && t<|\mathcal{T}|,\\
    \boldsymbol{\xi}_{s,t}
    &\geqslant \Phi_t(\cf)
      +c_t^{\text{penalty}}
       \left(P_t^{\min}(\cf)-d_{s,t}\right),\\
    \boldsymbol{\xi}_{s,t}
    &\geqslant \Phi_t(\cf)
      +c_t^{\text{penalty}}
       \left(R_t-H_t(\cf)\right),
\end{align}
\end{subequations}
for all applicable $s\in\mathcal{S}$ and $t\in\mathcal{T}$. The first two
inequalities follow by aggregating the startup and shutdown limits
\eqref{eq:UC:su-ramp-limits}--\eqref{eq:UC:sd-ramp-limits} with the power
balance and reserve requirement. The third aggregates committed minimum
output with the power balance, and the fourth aggregates reserve headroom with
\eqref{eq:UC:reserve-requirement}. Thus, any positive deficit on a right-hand
side requires at least that much penalized load imbalance or reserve
shortfall; when the deficit is nonpositive,
\eqref{eq:CRG:dispatch-floor-inequality} dominates. Adding the dispatch floor
$\Phi_t(\cf)$ therefore preserves every second-stage-feasible commitment while
strengthening the lower approximation of $Q_s(\x)$.

\paragraph{Hybrid valid-inequality generation.}
The price-based inequalities \eqref{eq:CRG:price-bound} are included in the
DW outer problem at initialization. After $K^{\mathrm{P}}$ DW pricing rounds,
the current outer solution is separated over
\eqref{eq:CRG:physical-recourse-bounds}; at most the $M^{\mathrm{R}}$ most
violated rows are added at each separation step and retained thereafter. 
This delay keeps the initial master compact and adds physical-recourse rows only when
the current commitment becomes physically inadequate.
Benders' cuts continue to be generated from the scenario inner problems.

Because the explicit commitment variables $\x$ remain in
\eqref{eq:UC:DW:general:master}, the dual effects of these rows pass to the
columns through the existing DW linking constraints. Hence the per-generator
pricing problem \eqref{eq:UC:DW:general:subproblem} requires no additional
variables or constraints.

\section{Computational Results}
\label{sec:computational-results}

The CRG algorithm is tested on a modified version of the FERC dataset from PGLIB-UC (\cite{PGLIB-UC2019,OnMIPFrmltnsForTheUCPrblm2018,RTOUCTstSystm2012}).
This dataset includes 934 thermal generators, aggregate wind generation, and 12 representative load scenarios.
Each load scenario covers 48 time steps.
Among the 934 thermal generators, 14 are CCs. 
In the original data set, these CCs are each modeled as single mode, traditional generators.
For this study, the data have been modified so that each CC is represented as a two combustion-turbine, one steam-turbine configuration with two operating modes.
Additional load scenarios are also generated using a copula-based method, as described in (\cite{WthrInfrmdPrblstcFrcstingAndScenGenInPS2024}).
The details of this load scenario generation are beyond the scope of this paper.
A penalty of \$50,000/MWh is applied to load shortage, excess generation, and reserve shortfall. Although actual penalty values vary by utility, such large values are standard in industry SUC models and contribute to the substantial difference observed between near-optimal and optimal decisions.

All tests are conducted on an Apple M3 Max CPU with 64 GB of RAM.
Gurobi version 12.0.2 is used as the solver for all inner, outer, and full problem formulations.
The CRG algorithm is compared against a BD algorithm, a Lagrangian Relaxation (LR) algorithm (see Appendix \cref{alg:LR}), and the B\&B approach on the full SUC problem.
Both CRG and BD use in-out separation with $\alpha = 0.4$ and $\beta=0.5$.
The BD and LR algorithms provide basic benchmarks for assessing the utility of the CRG algorithm for solving SUC with CCs.
It is important to note that no formulations or algorithms exist for solving SUC with CCs as of the publication of this paper.
All four algorithms are subject to a 3600-second time limit.
In Table \ref{sec:ComputationalResults:tab:time-and-solution}, the proposed algorithm (CRG) is italicized.

For the reported experiments with valid-inequalities within the CRG algorithm, 
we use $K^{\mathrm{P}}=17$ DW pricing rounds. At each separation step, at most
$M^{\mathrm{R}}=8$ physical rows are added.

\begin{figure}[]
    \centering
    \begin{subfigure}{.49\textwidth}
        \centering
        \includegraphics[width=\linewidth]{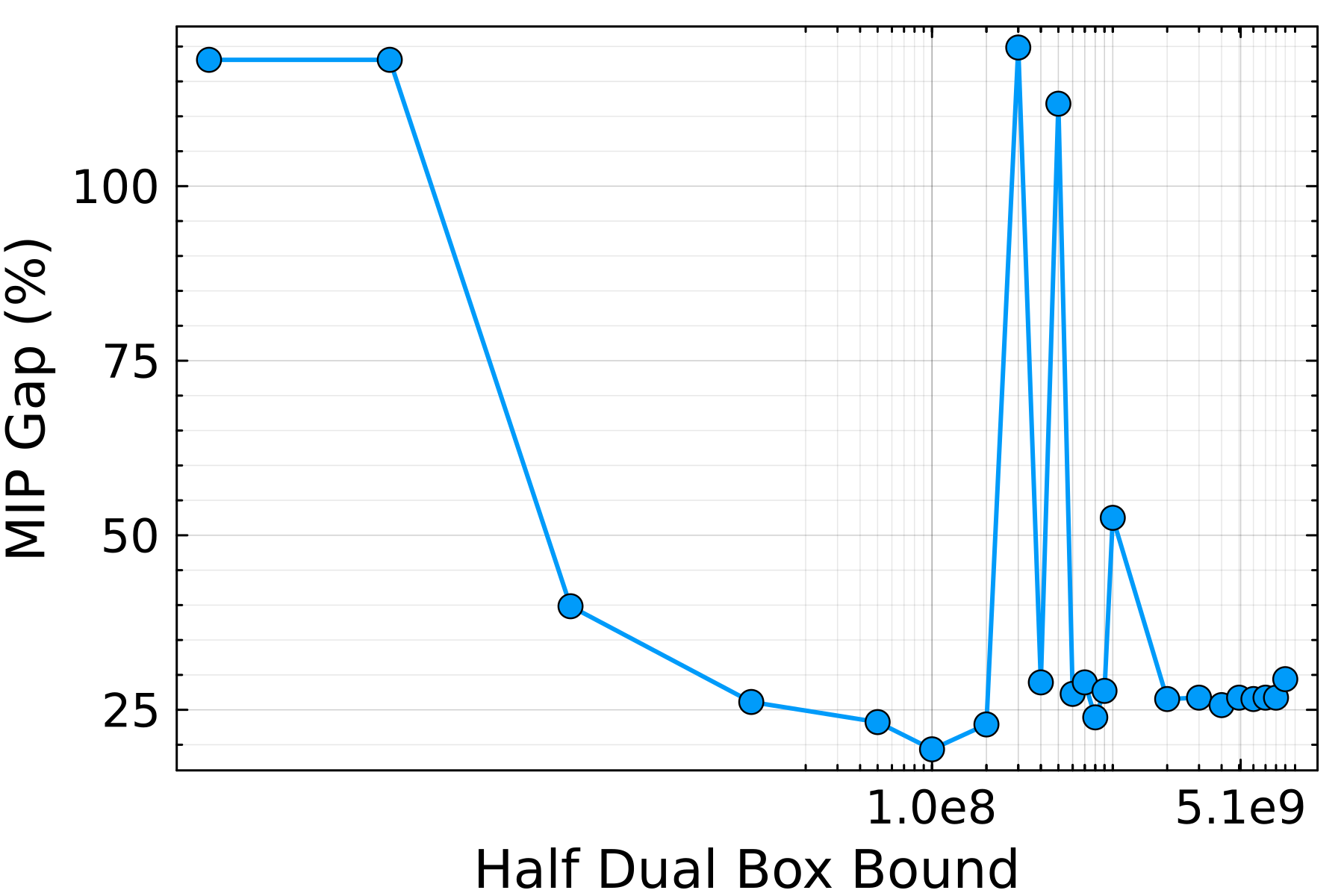}
        \caption{Dual-box half-width sensitivity.}
        \label{fig:dual-bound-sensitivity}
    \end{subfigure}
    \hfill
    \begin{subfigure}{.49\textwidth}
        \centering
        \includegraphics[width=\linewidth]{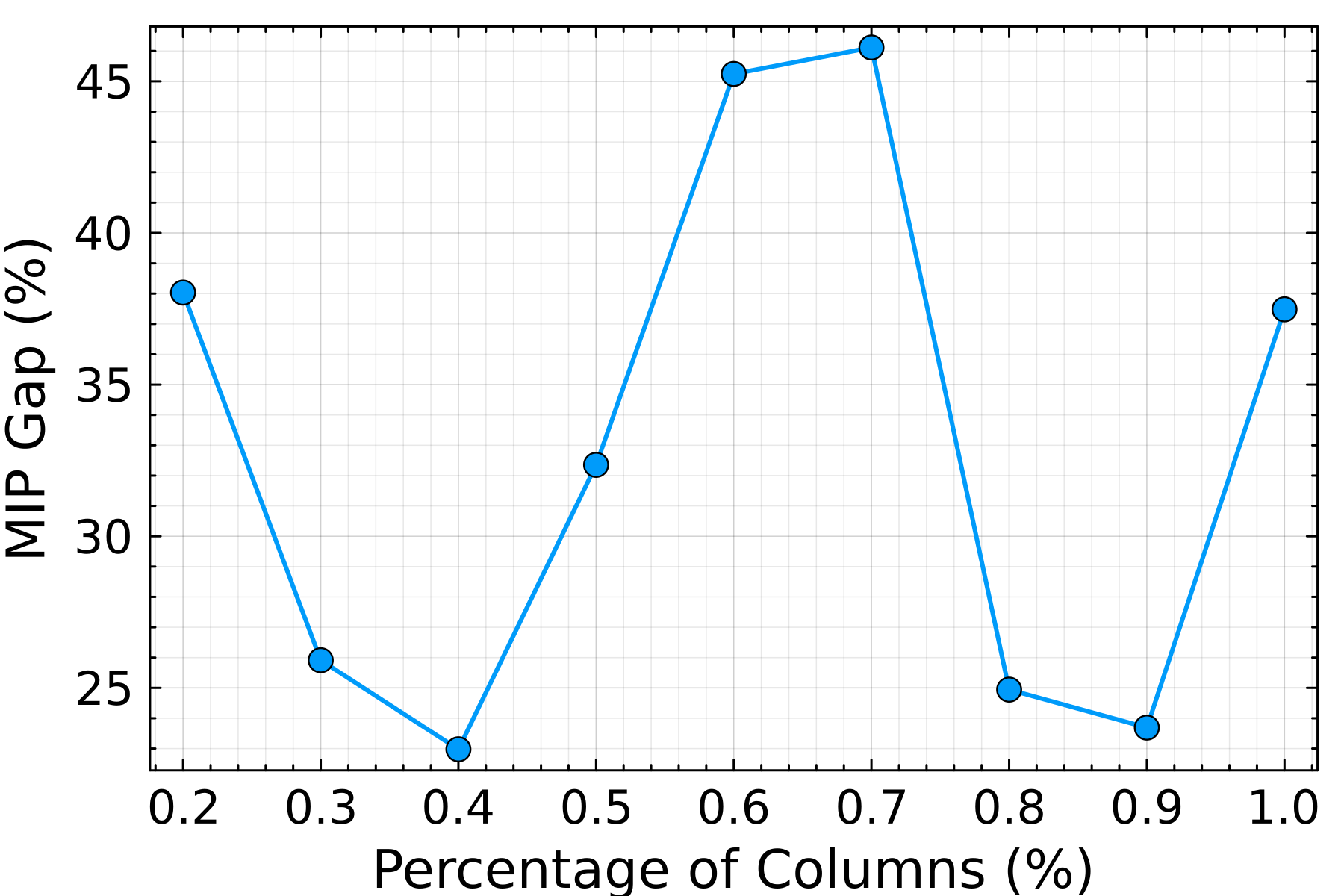}
        \caption{Column-addition percentage sensitivity.}
        \label{fig:col-perc-sensitivity}
    \end{subfigure}
    \caption{MIP gap sensitivity of CRG on the five-scenario instance to key algorithmic parameters: (a) the dual-stabilization box half-width, which limits deviations from the preceding dual solution, and (b) the percentage of columns with the smallest reduced costs selected for addition.}
    \label{fig:crg-parameter-sensitivity}
\end{figure}

Sensitivity analysis was performed on the same PGLIB-UC instance as the in-out separation sensitivity analysis from \cref{subsec:stabilization} to determine the best dual bound boxes for these instances (see \cref{fig:crg-parameter-sensitivity}(a)).
For these dual bound tests, the previously discovered best $\alpha$ and $\beta$ values were maintained constant for varying values of the box bounds.

The final tests were conducted using dual bound boxes of $\pm 1e8$ around the previous dual value.

Sensitivity analysis was also performed to determine what percentage of least-cost columns to add at each column generation step.
These tests were conducted on the PGLIB-UC instance with 5 scenarios (see \cref{fig:crg-parameter-sensitivity}(b)).

It is worth noting that these results also vary across the FERC dataset depending on the number of scenarios.
It was found through further testing that, despite the good performance of adding 40\% of the least cost columns for the 5 scenario case, adding 80\% of the least cost columns ended up having the best results across the higher scenario cases.
Thus, for the reported results later in this section, up to 80\% of the least cost columns were added at each column generation step of the CRG algorithm.

The CRG algorithm solves the outer problem faster than BD for 10, 30, and 40 scenarios, even with significant overhead (see Table \ref{sec:ComputationalResults:tab:time-and-solution}).
This performance improvement results from the significantly reduced search space generated by the pricing IPs.
Currently, the pricing IP schedule parsing takes a considerable amount of time, leading to an inflated total OP solve time. 
This total includes the time to parse and add the schedules to the OP. The time limit for each algorithm also accounts for the time required to parse and rewrite models within the algorithm. 
Consequently, the actual optimization solve time for CRG is approximately one tenth that of BD. 
With improved software engineering, this overhead could be significantly reduced, further enhancing the CRG algorithm’s ability to identify superior solutions within the same computational time frame.

\begin{table}[!t]
    \centering
    \caption{Performance comparison}
    \label{sec:ComputationalResults:tab:time-and-solution}
    \begin{tabular}{clrrrrr}
        \toprule
        &&\multicolumn{3}{c}{Solution Time (s)} & \\
        \cmidrule(lr){3-5}
        \#Scen. & Alg & \multicolumn{1}{c}{Outer} & \multicolumn{1}{c}{Inner} & Total & Obj. (\$bn) & Lb (\$bn)\\
        \midrule
        \phantom{0}
        10
            & \textbf{B\&B}
            & - & - & 3600
            & \textbf{5.17} & \textbf{5.16}\\
            & BD
            & 3213 & 132 & 3600
            & 5.62 & 4.88\\
            & LR
            & 3270 & 29 & 3600
            & 85.6 & 2.18\\
            & \textit{CRG}
            & \textit{3476} & \textit{124} & \textit{3600}
            & \textit{6.34} & \textit{4.54}\\
        \midrule
        20
            & \textbf{B\&B}
            & - & - & 3600
            & \textbf{5.18} & \textbf{5.16} \\
            & BD
            & 3115 & 255 & 3600
            & 7.33 & 4.75 \\
            & LR
            & 3239 & 47 & 3600
            & 90.7 & 2.19\\
            & \textit{CRG}
            & \textit{3532} & \textit{68} & \textit{3600}
            & \textit{6.92} & \textit{4.29} \\
        \midrule
        30  
            & \textbf{B\&B}
            & - & - & 3600
            & \textbf{5.42} & \textbf{5.17}\\
            & BD
            & 2902 & 424 & 3600
            & 9.37 & 4.69\\
            & LR
            & 3212 & 69 & 3600
            & 98.5 & 2.27\\
            & \textit{CRG}
            & \textit{3125} & \textit{475} & \textit{3600}
            & \textit{7.30} & \textit{3.78}\\
        \midrule
        40  & \textbf{B\&B}
            & - & - & 3600
            & \textbf{5.51} & \textbf{5.14} \\
            & BD
            & 2558 & 607 & 3600
            & 13.6 & 4.65\\
            & LR
            & 3212 & 108 & 3600
            & 106 & 2.26\\
            & \textit{CRG}
            & \textit{2912} & \textit{688} & \textit{3600}
            & \textit{10.0} & \textit{3.78} \\
        \midrule
        45  & B\&B
            & - & - & 3600
            & - & - \\
            & BD
            & 2770 & 715 & 3600
            & 18.6 & 4.57 \\
            & LR
            & 3324 & 134 & 3600
            & 111 & 2.30\\
            & \textbf{\textit{CRG}}
            & \textit{3442} & \textit{158} & \textit{3600}
            & \textbf{\textit{10.1}} & \textbf{\textit{3.78}} \\
        \midrule
        50 
            & B\&B
            & - & - & 3600
            & - & - \\
            & BD
            & 1412 & 802 & 3600
            & 18.8 & 4.56 \\
            & LR
            & 3132 & 150 & 3600
            & 113 & 2.33\\
            & \textbf{\textit{CRG}}
            & \textit{2489} & \textit{1111} & \textit{3600}
            & \textbf{\textit{10.0}} & \textbf{\textit{3.78}} \\
        \bottomrule
    \end{tabular}\\[0.5em]
    \footnotesize{Bolded instances indicate the algorithm with the smallest gap at the end of 3600s.}
\end{table}

\begin{figure}[!hbt]
\centering
\begin{subfigure}{.49\textwidth}
  \centering
  \includegraphics[width=\linewidth]{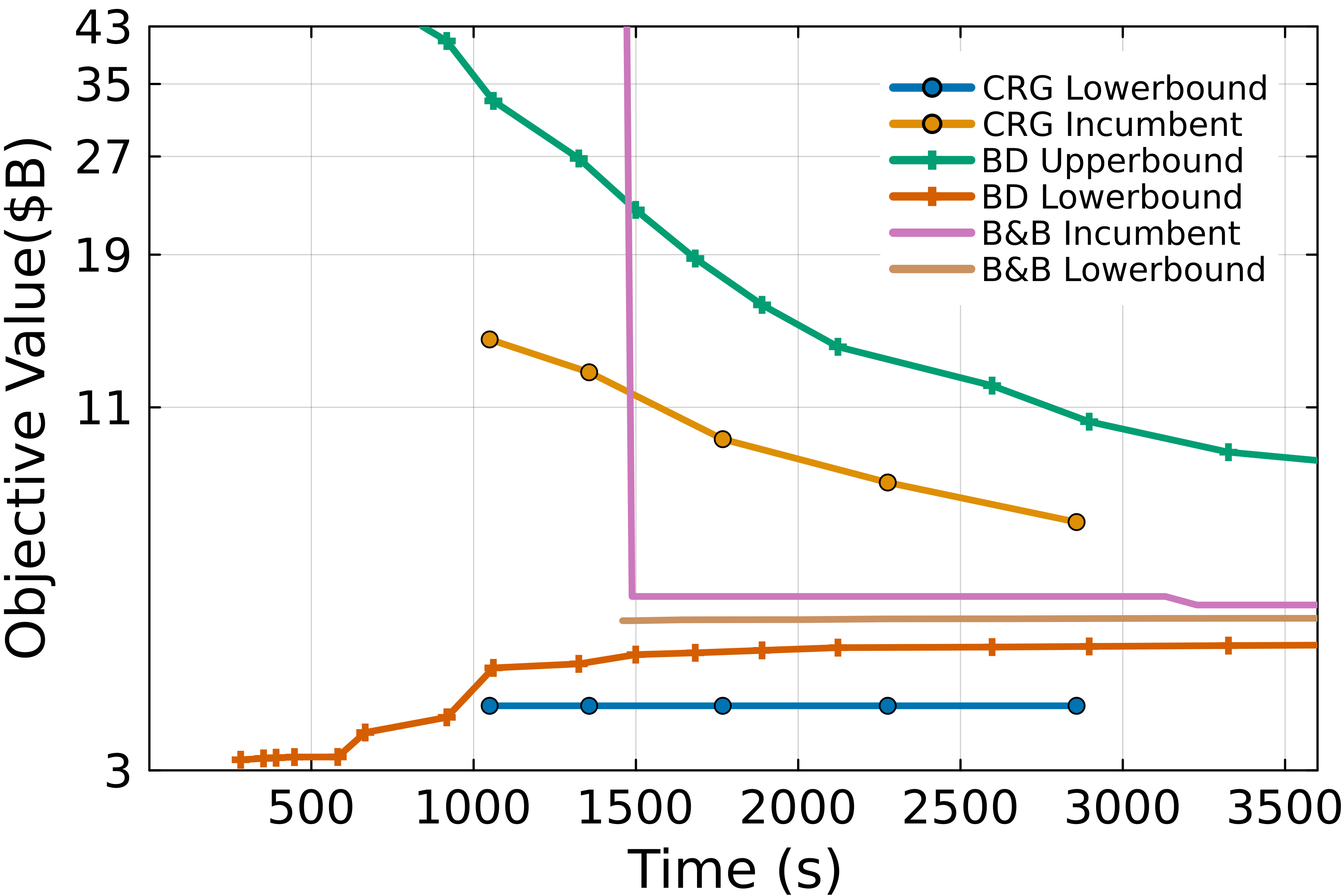}
  \caption{30 load scenarios}
  \label{fig:convergence-30-scenarios}
\end{subfigure}
\hfill
\begin{subfigure}{.49\textwidth}
  \centering
  \includegraphics[width=\linewidth]{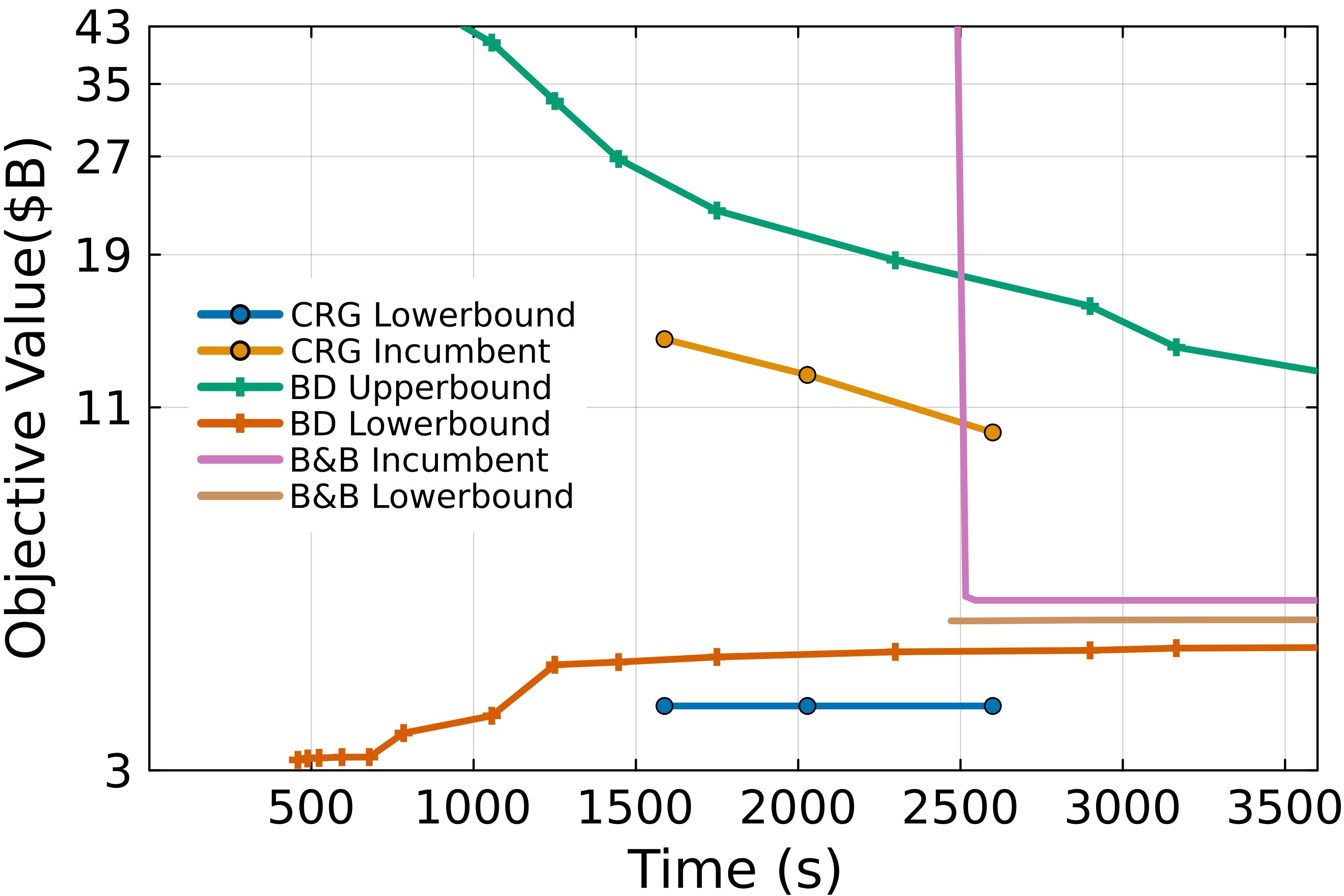}
  \caption{40 load scenarios}
  \label{fig:convergence-40-scenarios}
\end{subfigure} \\[1em]
\begin{subfigure}{.49\textwidth}
  \centering
  \includegraphics[width=\linewidth]{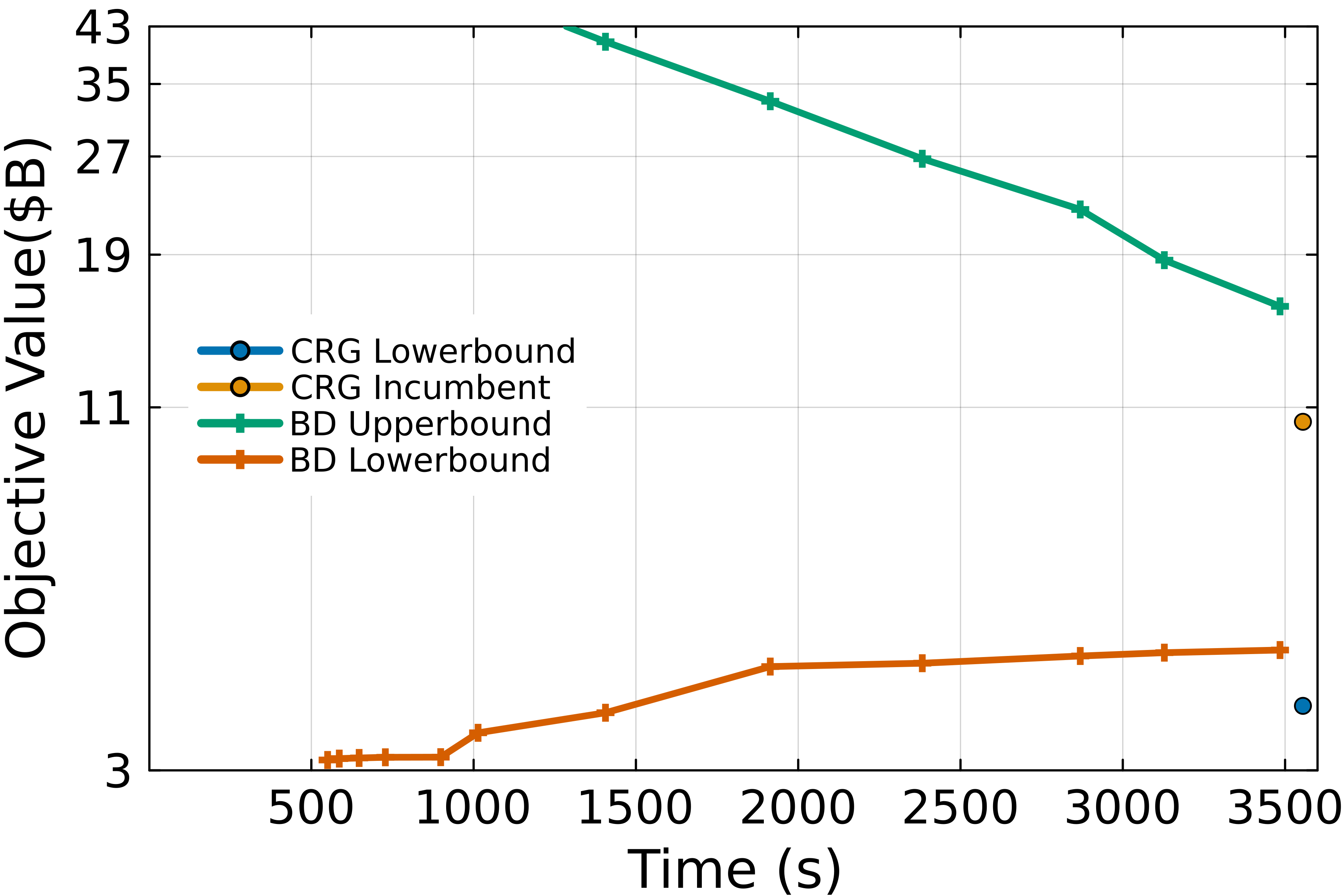}
  \caption{45 load scenarios}
  \label{fig:convergence-45-scenarios}
\end{subfigure}
\hfill
\begin{subfigure}{.49\textwidth}
  \centering
  \includegraphics[width=\linewidth]{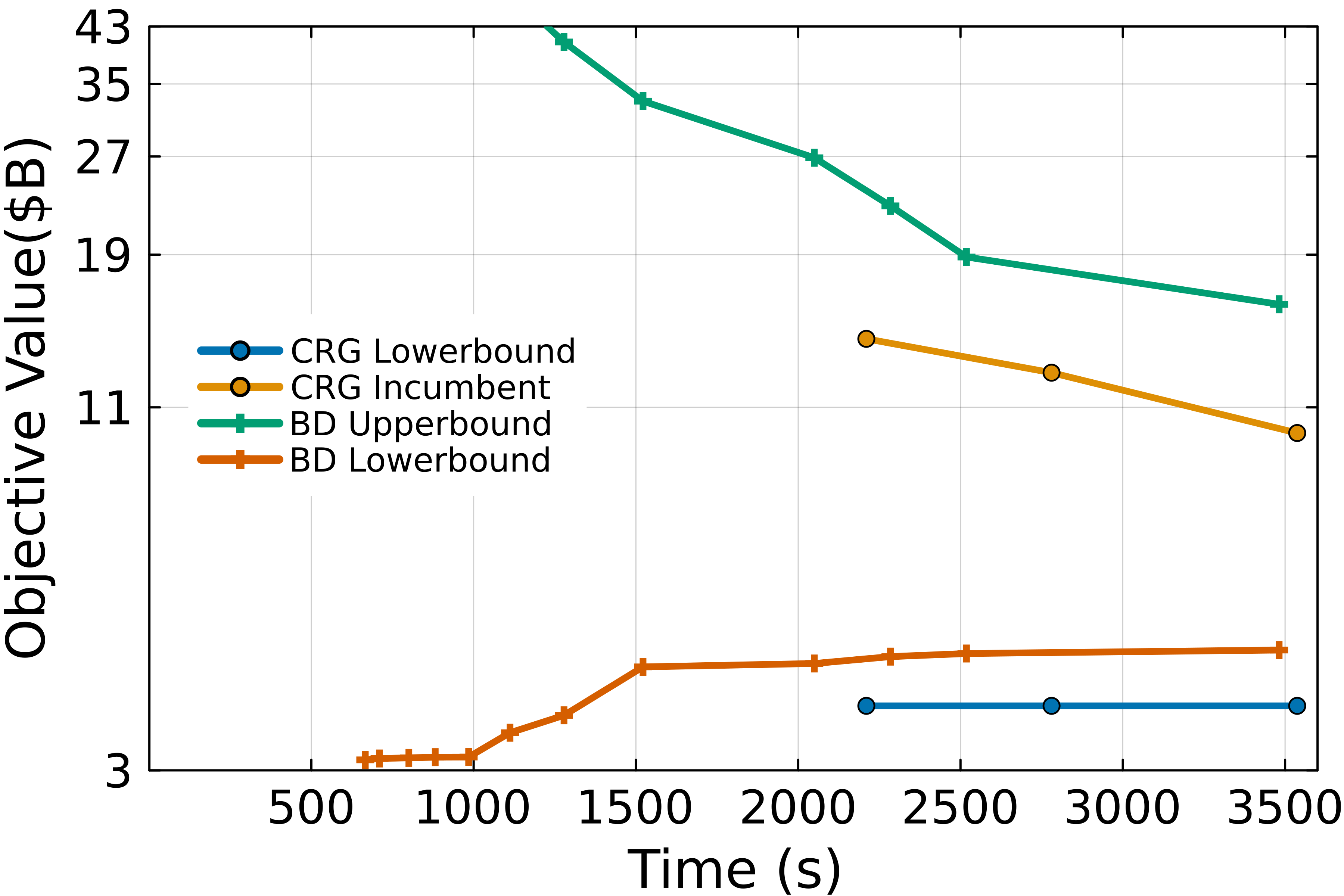}
  \caption{50 load scenarios}
  \label{fig:convergence-50-scenarios}
\end{subfigure}
\caption{Convergence of incumbent solutions and lower bounds for B\&B, BD, and CRG. Each plot shows a test case with a different number of load scenarios. For \cref{fig:convergence-45-scenarios} and \cref{fig:convergence-50-scenarios}, the B\&B solver ran out of memory before bounds could be found.}
\label{fig:convergence-plots}
\end{figure}

CRG finds better primal solutions faster than BD for all test cases above 10 scenarios, as seen by the final objective values in Table \ref{sec:ComputationalResults:tab:time-and-solution}.
For all load scenarios, CRG also provides tighter final lower bounds than LR.
The combination of high-quality primal solutions, tight lower bounds, and improved scalability leads to CRG's better convergence rates for larger test cases, as illustrated in \cref{fig:convergence-plots}.
These results demonstrate the efficacy of CRG in reducing the binary solution space to a smaller, higher-quality set of candidate solutions.

The valid inequalities introduced in \cref{subsec:hybrid-valid-inequalities} also greatly improve the lower bound of the CRG algorithm, as demonstrated during testing.
In initial tests of the CRG algorithm on 5 scenarios of the PG-LIB FERC dataset, CRG found a final lower bound of 1.078 \$B after one hour.
With valid inequalities, this lower bound was increased to 4.625 \$B for a one hour test on the same 5 scenarios.
This improvement demonstrates the efficacy of adding known cuts to improve convergence.
All results presented in this section benefit from the improved lower bounds achieved by these valid inequalities.

For instances with 45 or more scenarios, CRG produces tighter bounds than both B\&B and BD.
These instances are bolded in Table \ref{sec:ComputationalResults:tab:time-and-solution}.
At these problem sizes, even state-of-the-art B\&B solvers such as Gurobi suffer from the curse of dimensionality.
For example, at 40 scenarios, Gurobi must determine the optimal decision for 6,566,416 continuous and 128,128 integer variables after pre-solve.
At this size, Gurobi requires over 40GB of RAM to store the full set of variables, making this impractical for most ISOs, which often use standard work stations with 16 to 32GB of RAM. 
On such work stations, Gurobi would run out of memory after 30 scenarios.
In \cref{fig:convergence-45-scenarios} and \cref{fig:convergence-50-scenarios}, there is no incumbent value because of out-of-memory issues.

In summary, these observations highlight the importance of decomposition techniques such as CRG in scaling to larger system and scenario sizes. 
This is particularly relevant for SUC problems, where the large number of binary variables and constraints presents a major challenge for traditional B\&B solvers.

\section{Conclusion}
\label{sec:conclusion}

Motivated by the computational challenges of solving SUC problems with CCs, this paper has proposed a CRG algorithm that hybridizes Benders' and Dantzig-Wolfe decompositions.
This approach efficiently handles multiple scenarios through BD, while also tackling the combinatorial structure of CC operations through DW decomposition.
The addition of valid inequalities, in-out separation, and dual box bounds stabilizes and accelerates the convergence of CRG to improve the quality and speed of the algorithm over a baseline approach.
Numerical experiments demonstrate the efficiency and scalability of the proposed CRG algorithm.
Namely, compared to BD, CRG achieves better primal bounds, especially in the early stages of the optimization.
This behavior is particularly desirable for real-life settings, where high-quality solutions are sought within short computing times.
As the grid becomes larger and more volatile, this more scalable approach to solving SUC with CCs will be even more critical for daily grid operations.


The results reported in this paper use a modified FERC data set containing only 14 CCs --- 1.5\% of its generation mix --- whereas approximately 180 CCs would be needed to better represent the modern U.S. grid. Because CRG exploits generator separability, these experiments may understate its computational benefits on more representative instances. \textit{Future work} will therefore develop a more complete UC data set with mode-level CC information and a substantially larger share of CC generation, further refine the algorithmic framework, and improve software integration to enable large-scale SUC with CCs. In particular, we will investigate more efficient column-selection and valid-inequality-generation strategies, including machine-learning proxies, to reduce computational overhead.



%
\clearpage
\appendix
\section{Column-Row Generation Algorithm}
\begin{algorithm}[!ht]
\label{alg:CRG}
\caption{Column-Row Generation (CRG)}
\begin{algorithmic}[1]
\State \textbf{Input}: initial solution ($\x$), $\epsilon > 0$, $\alpha \in [0, 1]$, $\beta \in [0, 1]$
\State \texttt{corePnt} $\gets$ ($\x$), \texttt{ioPnt} $\gets$ ($\x$), \texttt{cutGenerated} $\gets$ \textbf{false}
\State $\bar{\Omega}$ $\gets$ $\{(\x)\}$, $\Gamma$ $\gets$ $\{\}$, UB $\gets$ $\infty$, LB $\gets$ $-\infty$
\While{(UB - LB)/UB $>$ $\epsilon$}
    \State \emph{$//$ Column Generation}
    \State Solve \eqref{eq:UC:DW:general:master}($\bar{\Omega}$, $\bar{\Gamma}$)
    \State Solve \eqref{eq:UC:DW:general:subproblem}, obtain column ($\x$)
    \State LB $\gets$ \eqref{eq:CRG:lower-bound}
    \If{ $\pi^{\top}\x + \sigma > 0$}
        \State $\bar{\Omega} \gets \bar{\Omega} \cup \{\x\}$
        \State Set \texttt{cutGenerated} $\gets$ \textbf{false}
        \State Continue
    \EndIf
    \State \emph{$//$ Row Generation}
    \If{\texttt{cutGenerated} is \textbf{false}}
        \State \texttt{ioPnt}$\gets$\eqref{eq:simple-io-separation}($\x$,\texttt{corePnt}, $\alpha$)
        \State Solve \eqref{eq:UC:BD:general:subproblem} and obtain cut ($\gamma$)
        \State $\bar{\Gamma} \gets \bar{\Gamma} \cup \{\gamma\}$
        \State \texttt{cutGenerated} $\gets$ \textbf{true}, \texttt{corePnt} $\gets$ \eqref{eq:corepoint}(\texttt{ioPnt}, $\beta$), UB $\gets$ \eqref{eq:CRG:upper-bound}($\omega, \y$)
    \Else
        \State \emph{$//$ Compute integer-feasible solution}
        \State Solve \eqref{eq:UC:DW:general:master}($\bar{\Omega}$, $\bar{\Gamma}$) with integer constraints
        \State Solve \eqref{eq:UC:BD:general:subproblem} and obtain cut ($Q$)
        \State add ($Q$) to $\bar{\Gamma}$
        \State \texttt{cutGenerated} $\gets$ \textbf{true}, \texttt{corePnt} $\gets$ \eqref{eq:corepoint}(\texttt{ioPnt}, $\beta$), UB $\gets$ \eqref{eq:CRG:upper-bound}($\omega, \y$)
    \EndIf
\EndWhile
\State \Return Optimal solution to \eqref{eq:UC:DW:general:master}
\end{algorithmic}
\end{algorithm}

\section{Lagrangian Relaxation Algorithm}
\begin{algorithm}[!ht]
\label{alg:LR}
\caption{Lagrangian Relaxation Algorithm (LR)}
\begin{algorithmic}[1]
\State \textbf{Input}: initial solution ($x_1^*$), initial step size $a$, $k = 1$
\State $\texttt{DualObjective}_1$ $\gets$ $\infty$, $\texttt{RelaxedObjective}_1$ $\gets$ $\infty$
\State $\x$ $\gets$ $x_1^*$ in Dual Problem (DP) 
\State Solve DP and obtain initial $\texttt{LagrangianMults}_1$
\State $\texttt{LRPenaltyTerms}_{1}$ $\gets$ $\texttt{LagrangianMults}_1$
\While{($\texttt{DualObjective}_k$ - $\texttt{RelaxedObjective}_k$)/$\texttt{DualObjective}_k$ $> \epsilon$}
    \State Solve LR and obtain solution $x_k^*$
    \State Update $\texttt{RelaxedObjective}_k$
    \State $\x$ $\gets$ $x_k^*$ in DP
    \State Solve DP and obtain solution $\texttt{LagrangianMults}_k$
    \State Update $\texttt{DualObjective}_k$
    \State $a_k$ $\gets$ $(1e-6)*\frac{\texttt{DualObjective}_k- \texttt{RelaxedObjective}_k}{\|\texttt{SubGradient}_k\|^2}$
    \State $\texttt{LRPenaltyTerms}_{k}$ $\gets$ $\texttt{LagrangianMults}_k + (a_k * \texttt{SubGradient}_k)$
    \State $k = k + 1$
\EndWhile
\end{algorithmic}
\end{algorithm}
%
%



\clearpage
\printbibliography

@ARTICLE{CompModeModforShortTermSchedofCCs2009,
  author={Liu, Cong and Shahidehpour, Mohammad and Li, Zuyi and Fotuhi-Firuzabad, Mahmoud},
  journal={IEEE Transactions on Power Systems}, 
  title={Component and Mode Models for the Short-Term Scheduling of Combined-Cycle Units}, 
  year={2009},
  volume={24},
  number={2},
  pages={976-990},
  doi={10.1109/TPWRS.2009.2016501}}

@article{StochUCProblem2004,
	author = {Shiina, Takayuki and Birge, John R.},
	doi = {https://doi.org/10.1111/j.1475-3995.2004.00437.x},
	eprint = {https://onlinelibrary.wiley.com/doi/pdf/10.1111/j.1475-3995.2004.00437.x},
	journal = {International Transactions in Operational Research},
	number = {1},
	pages = {19-32},
	title = {Stochastic unit commitment problem},
	url = {https://onlinelibrary.wiley.com/doi/abs/10.1111/j.1475-3995.2004.00437.x},
	volume = {11},
	year = {2004}
}

@misc{USElctrcGenByEngySrce2023,
    author = {U.S. Energy Information Administration},
    title = {What Is U.S. Electricity Generation by Energy Source},
    year = {2023},
    url = {https://www.eia.gov/tools/faqs/faq.php?id=427&t=3},
}

@techreport{ElctrcPwrAnnual2023,
    author = {Independent Statistics and Analysis},
    title = {Electric Power Annual 2023},
    institution = {U.S. Energy Information Administration},
    year = {2023},
    url = {https://www.eia.gov/electricity/annual/pdf/epa.pdf},
}

@INPROCEEDINGS{ASmplfdCCUnitMdlForMILPBasedUC2008,
  author={Chang, G. W. and Chuang, G. S. and Lu, T. K.},
  booktitle={2008 IEEE Power and Energy Society General Meeting - Conversion and Delivery of Electrical Energy in the 21st Century}, 
  title={A simplified combined-cycle unit model for mixed integer linear programming-based unit commitment}, 
  year={2008},
  volume={},
  number={},
  pages={1-6},
  doi={10.1109/PES.2008.4596496}}

@INPROCEEDINGS{AnEffcntApprchForUCAndEDwithCCandACPF2016,
  author={Bragin, Mikhail A. and Luh, Peter B. and Yan, Joseph H. and Stern, Gary A.},
  booktitle={2016 IEEE Power and Energy Society General Meeting (PESGM)}, 
  title={An efficient approach for Unit Commitment and Economic Dispatch with combined cycle units and AC Power Flow}, 
  year={2016},
  volume={},
  number={},
  pages={1-5},
  doi={10.1109/PESGM.2016.7741156}}

@ARTICLE{ANvlDcmpstinAndCoordntnAprchForLrgDAUCwithCCs2018,
  author={Sun, Xiaorong and Luh, Peter B. and Bragin, Mikhail A. and Chen, Yonghong and Wan, Jie and Wang, Fengyu},
  journal={IEEE Transactions on Power Systems}, 
  title={A Novel Decomposition and Coordination Approach for Large Day-Ahead Unit Commitment With Combined Cycle Units}, 
  year={2018},
  volume={33},
  number={5},
  pages={5297-5308},
  doi={10.1109/TPWRS.2018.2808272}}

@ARTICLE{ShrtTermSchdlngOfCCs2004,
  author={Bo Lu and Shahidehpour, M.},
  journal={IEEE Transactions on Power Systems}, 
  title={Short-term scheduling of combined cycle units}, 
  year={2004},
  volume={19},
  number={3},
  pages={1616-1625},
  doi={10.1109/TPWRS.2004.831706}}

@INPROCEEDINGS{SLRAndBCforUCwithCCs2014,
  author={Bragin, Mikhail A. and Luh, Peter B. and Yan, Joseph H. and Stern, Gary A.},
  booktitle={2014 IEEE PES General Meeting | Conference \& Exposition}, 
  title={Surrogate Lagrangian relaxation and branch-and-cut for unit commitment with combined cycle units}, 
  year={2014},
  volume={},
  number={},
  pages={1-5},
  doi={10.1109/PESGM.2014.6939901}}

@ARTICLE{PrceBsdUC:ACaseOfLRvsMIP2005,
  author={Tao Li and Shahidehpour, M.},
  journal={IEEE Transactions on Power Systems}, 
  title={Price-based unit commitment: a case of Lagrangian relaxation versus mixed integer programming}, 
  year={2005},
  volume={20},
  number={4},
  pages={2015-2025},
  doi={10.1109/TPWRS.2005.857391}}

@ARTICLE{UCwithFlxbleGenUnits2005,
  author={Bo Lu and Shahidehpour, M.},
  journal={IEEE Transactions on Power Systems}, 
  title={Unit commitment with flexible generating units}, 
  year={2005},
  volume={20},
  number={2},
  pages={1022-1034},
  doi={10.1109/TPWRS.2004.840411}}

@ARTICLE{HybridCmpntAndConfigModelForCCsInUC2018,
  author={Fang, Xin and Bai, Linquan and Li, Fangxing and Hodge, Bri-Mathias},
  journal={Journal of Modern Power Systems and Clean Energy}, 
  title={Hybrid component and configuration model for combined-cycle units in unit commitment problem}, 
  year={2018},
  volume={6},
  number={6},
  pages={1332-1337},
  doi={10.1007/s40565-018-0409-1}}

@article{BDandCRGforSlvngLrgSclLPsWithClDepRows2018,
	author = {{\.I}brahim Muter and {\c S}. {\.I}lker Birbil and Kerem B{\"u}lb{\"u}l},
	doi = {https://doi.org/10.1016/j.ejor.2017.06.044},
	issn = {0377-2217},
	journal = {European Journal of Operational Research},
	number = {1},
	pages = {29-45},
	title = {Benders decomposition and column-and-row generation for solving large-scale linear programs with column-dependent-rows},
	url = {https://www.sciencedirect.com/science/article/pii/S0377221717305854},
	volume = {264},
	year = {2018}}

@article{SlvngTwoStgeROPrblmsUsingCCGMthd2013,
	author = {Bo Zeng and Long Zhao},
	doi = {https://doi.org/10.1016/j.orl.2013.05.003},
	issn = {0167-6377},
	journal = {Operations Research Letters},
	number = {5},
	pages = {457-461},
	title = {Solving two-stage robust optimization problems using a column-and-constraint generation method},
	url = {https://www.sciencedirect.com/science/article/pii/S0167637713000618},
	volume = {41},
	year = {2013}}

@article{CmbngDWandBDtoSlveALrgSclNclrOutgePlnningPrblm2022,
	author = {Rodolphe Griset and Pascale Bendotti and Boris Detienne and Marc Porcheron and Halil {\c S}en and Fran{\c c}ois Vanderbeck},
	doi = {https://doi.org/10.1016/j.ejor.2021.07.018},
	issn = {0377-2217},
	journal = {European Journal of Operational Research},
	number = {3},
	pages = {1067-1083},
	title = {Combining Dantzig-Wolfe and Benders decompositions to solve a large-scale nuclear outage planning problem},
	url = {https://www.sciencedirect.com/science/article/pii/S0377221721006135},
	volume = {298},
	year = {2022}}

@article{CRGMthdToABtchMachneSchedPrblm2010,
author = {Wang, Gongshu and Tang, Lixin},
year = {2010},
month = {01},
pages = {},
title = {A Row-and-Column Generation Method to a Batch Machine Scheduling Problem},
journal = {Proceedings of the Ninth International Symposium on Operations Research and Its Applications (ISORA-10)}
}

@article{OnMIPFrmltnsForTheUCPrblm2018,
  author       = {Knueven, Bernard and Ostrowski, James and Watson, Jean-Paul},
  title        = {On Mixed Integer Programming Formulations for the Unit Commitment Problem},
  url          = {https://www.osti.gov/biblio/1492385},
  journal      = {Optimization Online Repository},
  issn         = {ISSN 9999-0042},
  volume       = {2018},
  place        = {United States},
  publisher    = {Mathematical Optimization Society},
  year         = {2018},
  month        = {11}}

@article{RTOUCTstSystm2012,
    author = {Krall, Eric and Higgins, Michael and O’Neill, Richard P.},
    title = {RTO Unit Commitment Test System},
    journal = {Federal Energy Regulatory Commission},
    year = {2012},
    url={https://www.ferc.gov/power-sales-and-markets/increasing-efficiency-through-improved-software/rto-unit-commitment-test}}

@Manual{PGLIB-UC2019,
    title = {PGLIB-UC},
    author = {Kneuven, Bernard and Coffrin, Carleton},
    year = {2019},
    note = {release v19.08},
    url = "https://github.com/power-grid-lib/pglib-uc",
  }

@misc{WthrInfrmdPrblstcFrcstingAndScenGenInPS2024,
      title={Weather-Informed Probabilistic Forecasting and Scenario Generation in Power Systems}, 
      author={Hanyu Zhang and Reza Zandehshahvar and Mathieu Tanneau and Pascal Van Hentenryck},
      year={2024},
      eprint={2409.07637},
      archivePrefix={arXiv},
      primaryClass={stat.ML},
      url={https://arxiv.org/abs/2409.07637}, 
}

@article{EfctvHybrdDecmpApprchToSlvNtwrkCnstrndStochUC2024,
title = {An effective hybrid decomposition approach to solve the network-constrained stochastic unit commitment problem in large-scale power systems},
journal = {EURO Journal on Computational Optimization},
volume = {12},
pages = {100085},
year = {2024},
issn = {2192-4406},
doi = {https://doi.org/10.1016/j.ejco.2024.100085},
url = {https://www.sciencedirect.com/science/article/pii/S2192440624000029},
author = {Ricardo M. Lima and Gonzalo E. Constante-Flores and Antonio J. Conejo and Omar M. Knio}
}

@inproceedings{ImplmntngAutoBDInAMdrnMIPSlvr2020,
	address = {Cham},
	author = {Bonami, Pierre and Salvagnin, Domenico and Tramontani, Andrea},
	booktitle = {Integer Programming and Combinatorial Optimization},
	editor = {Bienstock, Daniel and Zambelli, Giacomo},
	isbn = {978-3-030-45771-6},
	pages = {78--90},
	publisher = {Springer International Publishing},
	title = {Implementing Automatic Benders Decomposition in a Modern MIP Solver},
	year = {2020}}

@ARTICLE{AFstSltionMthdForLrgSclUCBsdOnLRAndDP2024,
  author={Hou, Jiangwei and Zhai, Qiaozhu and Zhou, Yuzhou and Guan, Xiaohong},
  journal={IEEE Transactions on Power Systems}, 
  title={A Fast Solution Method for Large-Scale Unit Commitment Based on Lagrangian Relaxation and Dynamic Programming}, 
  year={2024},
  volume={39},
  number={2},
  pages={3130-3140},
  doi={10.1109/TPWRS.2023.3287199}}

@ARTICLE{ARLEmbddSLRMthdForFstSlvngUCPrblms2025,
  author={Zhu, Yuhang and Cui, Gaochen and Liu, Anbang and Jia, Qing-Shan and Guan, Xiaohong and Zhai, Qiaozhu and Guo, Qi and Guo, Xianping},
  journal={IEEE Transactions on Power Systems}, 
  title={A Reinforcement Learning Embedded Surrogate Lagrangian Relaxation Method for Fast Solving Unit Commitment Problems}, 
  year={2025},
  volume={40},
  number={5},
  pages={3806-3818},
  doi={10.1109/TPWRS.2025.3529700}}

@book{BrnchAndPrce2026,
    author = {Desrosiers, Jacques and Lübbecke, Marco and Desaulniers, Guy and Gauthier, Jean Bertrand},
    title = {Branch-and-Price},
    publisher = {Springer Cham},
    year = {2026},
    pages={663},
    doi={https://doi.org/10.1007/978-3-031-96917-1},
    isbn={978-3-031-96917-1},
}





  



\end{document}